\makeatletter
\RequirePackage{keyval}
\define@key{Gm}{twosideshift}{%
	\ifdim#1=0pt%
	\else
	\PackageWarning{geometry}{Option `twosideshift` is gone since version 5}%
	\fi
}
\makeatother
\providecommand*{\StandardLayout}{}

\documentclass[envcountsame]{my-mmnp}

\usepackage{placeins}
\usepackage{caption}
\usepackage{mathtools}
\usepackage{amsfonts}
\usepackage{subfig}
\usepackage{pdfsync}
\usepackage{graphicx}
\usepackage{epstopdf}
\usepackage{morefloats}
\usepackage{amsmath,amssymb}
\usepackage{url}

\idline{ISSN: 0973-5348}{1}
\doi{\url{http://dx.doi.org/10.1051/mmnp/2018015}}

\numberwithin{equation}{section}

\begin{document}
	
\title{Uniform asymptotic stability\\ 
of a fractional tuberculosis model}

\author{Weronika Wojtak \inst{1, 2} 
\sep Cristiana J. Silva \inst{3}
\sep Delfim F. M. Torres \inst{3,}\thanks{\email{delfim@ua.pt}}} 


\vspace{0.5cm}

\institute{\inst{1}Algoritmi Center, University of Minho, 4800-058 Guimar\~{a}es, Portugal.\\
	\inst{2}Center of Mathematics, University of Minho, 4800-058 Guimar\~{a}es, Portugal.\\
	\inst{3}Center for Research and Development in Mathematics and Applications (CIDMA),\\
	Department of Mathematics, University of Aveiro, 3810--193 Aveiro, Portugal.}


\abstract{We propose a Caputo type fractional-order mathematical model 
for the transmission dynamics of tuberculosis (TB). Uniform asymptotic 
stability of the unique endemic equilibrium of the fractional-order TB 
model is proved, for any $\alpha \in (0, 1)$. Numerical simulations 
for the stability of the endemic equilibrium are provided.}






\titlerunning{Uniform asymptotic stability of a fractional tuberculosis model}

\authorrunning{W. Wojtak \sep C. J. Silva \sep D. F. M. Torres}

\maketitle


\section{Introduction}

When defining a differential operator, one often employs, besides ordinary derivatives, 
generalized derivatives, which appear in a natural way when considering extensions 
of differential operators defined on differentiable functions, and weak derivatives, 
related to the transition to the adjoint operator \cite{MyID:298}. Derivatives of fractional 
and negative orders appear when the differentiation is defined by means 
of an integral transform, applicable to the domain of definition and range 
of such generalized differential operator \cite{MyID:310}. 
This is often done in order to obtain the simplest possible representation 
of the corresponding differential operator of a function and to attain 
a reasonable generality in the formulation of problems 
and satisfactory properties of the objects considered \cite{book:frac}. 
Problems in the theory of differential equations, e.g., problems of existence, 
uniqueness, regularity, continuous dependence of the solutions 
on the initial data or on the right-hand side, or the explicit form 
of a solution of a differential equation defined by a given differential expression, 
are readily interpreted in the theory of operators as problems on the corresponding 
differential operator defined on suitable function spaces \cite{MR3200762}. 
One advantage of fractional order differential equations is that they provide 
a powerful instrument to incorporate memory and hereditary properties 
into the systems, as opposed to the integer order models, where such effects 
are neglected or difficult to incorporate \cite{GonzalezParra}. Moreover, 
in order to precisely reproduce the nonlocal, frequency- and history-dependent 
properties of power law phenomena, some different modeling tools, 
based on fractional operators, have to be introduced: see, e.g., 
\cite{Atangana} and references therein.

Communicable diseases have always been an important part of human history \cite{Brauer:17}.
Fractional-order differential system models for infectious disease dynamics have been 
introduced in recent years \cite{AhmedPhysicaA2007,DingYeFracHIVCD4cells,GonzalezParra,%
FracCarlaPintoMalaria,SofiaShakoorDelfim,ZhangDDNS2015}.
In \cite{SofiaShakoorDelfim}, the authors propose a fractional-order model and show, 
through numerical simulations, that the fractional models fit better the first dengue 
epidemic recorded in the Cape Verde islands off the coast of west Africa, in 2009, 
when compared with the standard differential model. The authors of \cite{GonzalezParra} 
show that a nonlinear fractional order epidemic model is well suited to provide numerical 
results that agree very well with real data of influenza A (H1N1) at the population level. 
In \cite{FracCarlaPintoMalaria}, a fractional model for malaria transmission is considered 
and numerical simulations are done for the variation of the values of the fractional 
derivative and of the parameter that models personal protection. In \cite{ZhangDDNS2015}, 
fractional-order derivatives are introduced into an HIV infection model and local 
asymptotic stability is proved. The authors of \cite{DingYeFracHIVCD4cells} introduce 
fractional-order derivatives into a model of HIV infection of $CD4^+$ T-cells  
and analyze the local asymptotic stability of the equilibrium points. 
In \cite{VargasDeLeon201575}, the uniform asymptotic stability is proved 
extending the Volterra-type Lyapunov functions to fractional-order epidemic systems. 
Fractional-order predator-prey models are investigated in \cite{AhmedMMA2007}.
In particular, existence and uniqueness of solutions are proved, 
and stability of equilibrium points studied. 
Numerical solutions of such models were obtained \cite{AhmedMMA2007}. 
Here we are interested to investigate fractional calculus
with respect to tuberculosis.

Tuberculosis (TB) is a bacterial disease caused by \emph{Mycobacterium tuberculosis}, 
which is usually spread through the air by people with active TB. 
TB is one of the top ten causes of death worldwide, which justifies
the amount of research on the area: see, e.g., \cite{MyID:353,MyID:271,reportWHO2016}.
For a good survey on optimal models of TB, we refer the reader to \cite{MyID:301}.
In \cite{MyID:353}, delays are introduced in a TB model, 
representing the time delay on the diagnosis and commencement 
of treatment of individuals with active TB infection. 
Optimal control strategies to minimize the cost of interventions, 
considering reinfection and post-exposure interventions, are investigated in \cite{MyID:271}.
In \cite{MyID:305}, the potential of two post-exposure interventions, 
treatment of early latent TB individuals and prophylactic treatment/vaccination 
of persistent latent TB individuals, is investigated.
A mathematical model for TB is studied in \cite{MyID:314}
from the optimal control point of view, using a multiobjective approach.
For numerical simulations using TB data from Angola see \cite{MyID:230}.
Despite the numerous works on tuberculosis, the literature on fractional-order 
mathematical models for TB is scarce.  In \cite{Sweilama:2016}, 
the authors propose a multi-strain TB model 
of variable-order fractional derivatives and develop a numerical scheme 
to approximate the endemic solution numerically. 
Here, we propose a Caputo type fractional-order mathematical model 
for the transmission dynamics of tuberculosis (TB), based on the nonlinear 
differential system studied in \cite{Yang201079}, and investigate its stability. 

The stability question is of main interest in epidemiological 
control systems \cite{Matignon1996,MyID:355,MyID:373}.
For nonlinear fractional-order systems, stability analysis is a recent and interesting topic: see, e.g., \cite{Baleanu:2010,Chen2014633,Delavari2012,LiChen:Automatica:2009,LiChen:CMA:2010,LiMa2013,%
LiZhang:2011,Rivero:2013,VargasDeLeon201575}. The direct Lyapunov method provides a way 
to determine (asymptotic, uniform, global asymptotic) stability of an equilibrium point 
without explicitly solving for integer-order nonlinear systems \cite{VargasDeLeon201575}.  
In \cite{Delavari2012}, the authors extended the Lyapunov direct method to Caputo type 
fractional-order nonlinear systems using Bihari's inequality and Bellman--Grownwall's inequality. 
Different approaches for the  extension of the Lyapunov direct method were proposed 
in \cite{Baleanu:2010,LiChen:Automatica:2009,LiChen:CMA:2010}. 
See \cite{LiZhang:2011,Petras:2009,Petras:2011,Rivero:2013,Trigeassou:2011} 
for a survey on stability analysis of fractional differential equations.
In \cite{Sweilama:2016}, the stability of the endemic equilibrium of a 
fractional-order TB model is studied numerically. Here we prove, analytically, 
the uniform asymptotic stability of the unique endemic equilibrium 
of the fractional-order TB model, for any $\alpha \in (0, 1)$. Moreover, 
we illustrate the theoretical stability results through numerical simulations.  

The paper is organized as follows. In Section~\ref{sec:pre}, we present basic 
definitions and some known results about Caputo fractional calculus 
as well as results on uniform asymptotic stability and Volterra-type 
Lyapunov functions for fractional-order systems. In Section~\ref{sec:frac:model}, 
we introduce Caputo fractional-order derivatives into the TB model proposed 
in \cite{Yang201079} and study the existence of equilibrium points. 
In Section~\ref{sec:unif:stab}, we prove the uniform asymptotic stability 
of the unique endemic equilibrium of the fractional-order TB model. 
Numerical simulations are provided in Section~\ref{sec:numsim}, 
which illustrate the stability result proved in the previous section. 
We finish with Section~\ref{sec:conc} of conclusions. 


\section{Preliminaries on the Caputo fractional calculus}
\label{sec:pre}

We begin by introducing the definition of Caputo fractional derivative
and recalling its main properties.

\begin{dfntn}[See \cite{GJI:GJI529}]
Let $a > 0$, $t > a$, $\alpha, a, t \in\mathbb{R} $. The Caputo fractional 
derivative of order $\alpha$ of function $f \in C^n$ is given by
\begin{equation}
\label{Caputo}
_{a}^{C}D_{t}^{\alpha}f(t)= \dfrac{1}{\Gamma(n-\alpha)} 
\int_{a}^{t}\dfrac{f^{(n)}(\xi)}{(t-\xi)^{\alpha+1-n}}d\xi,
\qquad n-1<\alpha<n \in \mathbb{N}.
\end{equation}
\end{dfntn}

\begin{prprty}[Linearity, see, e.g., \cite{Diethelm}] 
Let $f(t),g(t):[a,b] \rightarrow \mathbb{R}$ be such that 
$_{a}^{C}D_{t}^{\alpha}f(t)$ and $_{a}^{C}D_{t}^{\alpha}g(t)$ 
exist almost everywhere and let $c_1,c_2 \in \mathbb{R}$. Then, 
$_{a}^{C}D_{t}^{\alpha}(c_1 f(t)+ c_2 g(t))$ exists almost 
everywhere, and
\begin{equation}
_{a}^{C}D_{t}^{\alpha} (c_{1} f(t)+ c_{2} g(t))
= c_{1} \vspace{1mm} _{a}^{C}D_{t}^{\alpha}f(t) 
+ c_{2}\vspace{1mm} _{a}^{C}D_{t}^{\alpha}g(t).
\end{equation}
\end{prprty}

\begin{prprty}[Caputo derivative of a constant, see, e.g., \cite{Podlubny}] 
The fractional derivative of a constant function $f(t)\equiv c$ is zero:
\begin{equation}
_{t_{0}}^{C}D_{t}^{\alpha} c = 0.
\end{equation}
\end{prprty}

Let us consider the following general fractional differential 
equation involving the Caputo derivative:
\begin{equation}
\label{CaputoGeneral}
_{a}^{C}D_{t}^{\alpha} x(t)= f(t,x(t)), 
\qquad \alpha \in (0,1),
\end{equation}
with the initial condition $x_0=x(t_0)$.

\begin{dfntn}[See, e.g., \cite{LiChen:Automatica:2009}]
\label{def:eq}
The constant $x^*$ is an equilibrium point of the Caputo fractional 
dynamic system \eqref{CaputoGeneral} if, and only if, $f(t, x^*) = 0$. 
\end{dfntn}

Next theorem is an extension of the Lyapunov direct method for 
Caputo type fractional-order nonlinear systems \cite{Delavari2012}. 

\begin{thrm}[Uniform Asymptotic Stability \cite{Delavari2012}]
\label{uniform_stability}
Let $x^*$ be an equilibrium point for the nonautonomous fractional 
order system \eqref{CaputoGeneral} and $\Omega \subset \mathbb{R}^{n}$ 
be a domain containing $x^*$. Let $L:[0,\infty) \times \Omega 
\rightarrow \mathbb{R}$ be a continuously differentiable function
such that
\begin{equation}
W_1(x) \leq L(t,x(t)) \leq W_2(x)
\end{equation}
and
\begin{equation}
_{a}^{C}D_{t}^{\alpha}  L(t,x(t)) \leq -W_3(x)
\end{equation}
for all $\alpha\in (0,1)$ and all $x \in \Omega$, 
where $W_1(\cdot)$, $W_2(\cdot)$ and $W_3(\cdot)$ 
are continuous positive definite functions on $\Omega$. 
Then the equilibrium point of system \eqref{CaputoGeneral} 
is uniformly asymptotically stable.
\end{thrm} 

In what follows we recall a lemma proved in \cite{VargasDeLeon201575}, 
where a Volterra-type Lyapunov function is obtained 
for fractional-order epidemic systems. 

\begin{lmm}[See \cite{VargasDeLeon201575}]
\label{lemma31}
Let $x(\cdot)$ be a continuous and differentiable function
with $x(t)\in \mathbb{R_{+}}$. Then, for any time instant 
$t \geq t_0$, one has
\begin{equation}
\label{lemma31eq}
_{t_{0}}^{C}D_{t}^{\alpha}\left[ x(t)-x^{*}-x^{*}
\ln\dfrac{x(t)}{x^{*}}\right]
\leq \left( 1-\dfrac{x^{*}}{x(t)}\right)
_{t_{0}}^{C}D_{t}^{\alpha} x(t), 
\qquad x^{*} \in \mathbb{R}^{+}, 
\qquad \forall\alpha\in (0,1).
\end{equation}
\end{lmm}


\section{Fractional-order tuberculosis (TB) model}
\label{sec:frac:model}

In this section we propose a Caputo fractional-order version 
of a tuberculosis (TB) model in \cite{Yang201079}. 
The model describes the dynamics of a population that is susceptible 
to infection by the \emph{Mycobacterium tuberculosis} with 
incomplete treatment. The total population is partitioned 
into four compartments: 
\begin{itemize}

\item susceptible individuals, $S$; 

\item latent individuals, $L$, which have been infected but are not infectious 
and do not have symptoms of the disease; 

\item infectious individuals, $I$, which have active TB, 
may transmit the infection but are not in treatment; 

\item and under treatment infected individuals, $T$.
\end{itemize}
The susceptible population is increased by the recruitment of individuals 
into the population, assumed susceptible, at a rate $\Lambda$. 
All individuals suffer from natural death, at a constant rate $\mu$. 
Susceptible individuals acquire TB infection by the contact 
with individuals in the class $I$ at a rate $\beta I S$, 
where $\beta$ is the transmission coefficient. Individuals 
in the latent class $L$ become infectious at a rate $\epsilon$, 
and infectious individuals $I$ start treatment at a rate $\gamma$. 
Treated individuals $T$ leave their compartment at a rate $\delta$. 
After leaving the treatment compartment, an individual may enter compartment 
$L$, due to the remainder of \emph{Mycobacterium tuberculosis}, 
or compartment $I$, due to the failure of treatment. The parameter $k$,
$0\leq k\leq 1$, represents the failure of the treatment, where $k=0$ means that all 
the treated individuals will become latent, while $k=1$ 
means that the treatment fails and all the treated individuals 
will still be infectious. Infectious, $I$, and under treatment individuals, $T$, 
may suffer TB-induced death at the rates $\alpha_1$ and $\alpha_2$, respectively. 
The Caputo fractional-order system that describes 
the previous assumptions is 
\begin{equation}
\label{model_TB_fractional}
\begin{cases}
_{t_{0}}^{C}D_{t}^{\alpha}S(t) = \Lambda-\beta I(t) S(t) - \mu S(t),\\
_{t_{0}}^{C}D_{t}^{\alpha}L(t) = \beta I(t) S(t) + (1-k)\delta T(t) - (\mu + \epsilon) L(t),\\
_{t_{0}}^{C}D_{t}^{\alpha}I(t) = \epsilon L(t) + k \delta T(t) - (\mu + \gamma + \alpha_{1}) I(t),\\
_{t_{0}}^{C}D_{t}^{\alpha}T(t) = \gamma I(t) - (\mu + \delta + \alpha_{2}) T(t).
\end{cases}
\end{equation}
The feasible region of system \eqref{model_TB_fractional} is given by
\begin{equation}
\label{omega}
\Omega = \left\{ (S,L,I,T) \in \mathbb{R}^{4}_{+}: S+L+I+T 
\leq  \dfrac{\Lambda}{\mu} \right\}.
\end{equation}
Let $b_{1}=\mu+\epsilon$, $b_{2}=\mu+\gamma + \alpha_{1}$ 
and $b_{3}=\mu+\delta + \alpha_{2}$. Consider the following 
endemic threshold \cite{Yang201079}: 
\begin{equation}
\label{R0}
R_{0}=\dfrac{\beta\epsilon b_{3}\Lambda}{\mu b_{1}b_{2}b_{3} 
- \mu \delta \gamma ((1-k)\epsilon + k b_{1})}.
\end{equation}
According to Definition~\ref{def:eq}, system \eqref{model_TB_fractional} 
has a disease-free equilibrium 
$$
E_{0} = \left( \dfrac{\Lambda}{\mu},0,0,0\right)
$$ 
and, when $R_{0} > 1$, an unique endemic equilibrium 
$$
E^{*}=(S^{*},L^{*},I^{*},T^{*})
$$ 
with
\begin{equation}
\label{EE}
S^{*} = \dfrac{\Lambda}{\mu + \beta I^{*}},\quad
L^{*} = \dfrac{1}{\epsilon}\left(b_{2}-\dfrac{k\delta\gamma}{b_{3}}\right)I^{*},\quad
I^{*} = \dfrac{\mu}{\beta}\left( R_{0}-1 \right) ,\quad
T^{*} = \dfrac{\gamma}{b_{3}}I^{*}.
\end{equation}


\section{Uniform asymptotic stability of the endemic equilibrium}
\label{sec:unif:stab}

In this section we prove uniform asymptotic stability of the endemic equilibrium  
$E^{*}=(S^{*},L^{*},I^{*},T^{*})$ of the fractional order system 
\eqref{model_TB_fractional}.

\begin{thrm}
\label{thm:MR}
Let $\alpha \in (0,1)$ and $R_0 > 1$.
Then the unique endemic equilibrium $E^*$ of the fractional order 
system \eqref{model_TB_fractional} is uniformly asymptotically stable 
in the interior of $\Omega$ defined by \eqref{omega}.
\end{thrm}

\begin{proof}
Consider the following Lyapunov function:
$$
V(t) = m_1 V_1(S(t)) + m_1 V_2(L(t)) + m_2 V_3(I(t) ) + m_3 V_4(T(t)),
$$
where
\begin{equation*}
\begin{split}
V_{1}(S(t))&= S(t) - S^* -S^*  \ln \dfrac{S(t)}{S^*},\\
V_{2}(L(t))&=L(t) - L^* - L^*\ln \dfrac{L(t)}{L^*},\\
V_{3}(I(t))&= I(t) - I^* - I^*\ln \dfrac{I(t)}{I^*},\\  
V_{4}(T(t))&= T(t) - T^* - T^* \ln \dfrac{T(t)}{T^*}
\end{split}
\end{equation*}
and 
$$
m_1=\epsilon, \quad
m_2=b_1,\quad 
m_3=\dfrac{\delta T^* (b_1 k + \epsilon (1-k))}{\gamma I^*}.
$$
Function $V$ is defined, continuous and positive definite 
for all $S(t)>0$, $L(t)>0$, $I(t)>0$ and $T(t) > 0$.
By Lemma~\ref{lemma31}, we have
\begin{equation*}
_{t_{0}}^{C}D_{t}^{\alpha}V \leq m_{1}\left( 1-\dfrac{S^{*}}{S}\right)
{_{t_{0}}^{C}D}_{t}^{\alpha} S + m_{1}\left( 1-\dfrac{L^{*}}{L}\right) 
{_{t_{0}}^{C}D}_{t}^{\alpha} L
+ m_{2}\left( 1-\dfrac{I^{*}}{I}\right)\vspace{0.5mm} _{t_{0}}^{C}D_{t}^{\alpha} I 
+ m_{3}\left( 1-\dfrac{T^{*}}{T}\right) {_{t_{0}}^{C}D}_{t}^{\alpha} T. 
\end{equation*}
It follows from \eqref{model_TB_fractional} that
\begin{equation*}
\begin{split}
_{t_{0}}^{C}D_{t}^{\alpha}V
&\leq m_{1}\left( 1-\dfrac{S^{*}}{S}\right) (\Lambda-\beta I S - \mu S) 
+ m_{1}\left( 1-\dfrac{L^{*}}{L}\right) (\beta I S + (1-k)\delta T - b_{1} L)\\
&\quad + m_2 \left( 1-\dfrac{I^*}{I}\right) \left( \epsilon L + k \delta T - b_2 I\right) 
+ m_3 \left(1-\dfrac{T^*}{T} \right)\left( \gamma I - b_3 T\right). 
\end{split}
\end{equation*}
At the endemic equilibrium $\Lambda = S^* \left( \mu + \beta I\right)$, we have
\begin{equation*}
\begin{split}
_{t_{0}}^{C}D_{t}^{\alpha}V
&\leq m_{1} \left[ \mu S^* \left(2 - \frac{S}{S^*} 
- \frac{S^*}{S} \right) + \beta I^* S^* \left( 1- \frac{S^*}{S} 
- \frac{S}{S^*}\frac{I}{I^*} + \frac{I}{I^*}  \right) \right]\\ 
&\quad + m_{1}\left( 1-\dfrac{L^{*}}{L}\right) \left(\beta I S + (1-k)\delta T - b_{1} L\right)\\
&\quad + m_2 \left( 1-\dfrac{I^*}{I}\right) \left( \epsilon L + k \delta T - b_2 I\right) 
+ m_3 \left(1-\dfrac{T^*}{T} \right)\left( \gamma I - b_3 T\right). 
\end{split}
\end{equation*}
Let $m_1 = \varepsilon$. After some simplifications, one has 
\begin{equation*}
\begin{split}
_{t_{0}}^{C}D_{t}^{\alpha}V
&\leq \varepsilon  \mu S^* \left(2 - \frac{S}{S^*} - \frac{S^*}{S} \right) 
+ \varepsilon \beta I^* S^* \left( 1- \frac{S^*}{S} + \frac{I}{I^*} 
-  \dfrac{L^{*}}{L} \frac{I}{I^*} \frac{S}{S^*} \right)\\ 
&\quad + \varepsilon (1-k)\delta T - \varepsilon b_{1} L
- \varepsilon\dfrac{L^{*}}{L} (1-k)\delta T + \varepsilon L^{*} b_{1}\\
&\quad + m_2 \left( 1-\dfrac{I^*}{I}\right) \left( \epsilon L + k \delta T 
- b_2 I\right) + m_3 \left(1-\dfrac{T^*}{T} \right)\left( \gamma I - b_3 T\right). 
\end{split}
\end{equation*}
Let $m_2 = b_1$. Using 
$$
b_1 L^* = \beta I^* S^* + (1-k) \delta T^*
$$ 
at the endemic equilibrium, we have
\begin{equation*}
\begin{split}
_{t_{0}}^{C}D_{t}^{\alpha}V
&\leq \varepsilon  \mu S^* \left(2 - \frac{S}{S^*} - \frac{S^*}{S} \right) 
+ \varepsilon \beta I^* S^* \left( 2- \frac{S^*}{S} + \frac{I}{I^*} 
- \frac{L}{L^*} -  \dfrac{L^{*}}{L} \frac{I}{I^*} \frac{S}{S^*} \right)\\ 
&\quad + \varepsilon (1-k)\delta T^* \left( 1+ \frac{T}{T^*} - \frac{L}{L^*} 
- \dfrac{L^{*}}{L} \frac{T}{T^*} \right)
+ m_3 \left(1-\dfrac{T^*}{T} \right)\left( \gamma I - b_3 T\right)
\epsilon b_1 L^* \frac{L}{L^*}\\
&\quad + b_1 k \delta T - b_1 b_2 I - \epsilon b_1 
L^*\dfrac{I^*}{I}  \frac{L}{L^*} - b_1\dfrac{I^*}{I} k \delta T + b_1 I^* b_2.
\end{split}
\end{equation*}
Moreover, since 
$$
b_1 L^* = \beta I^* S^* + (1-k) \delta T^*
$$ 
at the equilibrium point, it follows that 
\begin{equation*}
\begin{split}
_{t_{0}}^{C}D_{t}^{\alpha}V
&\leq \varepsilon  \mu S^* \left(2 - \frac{S}{S^*} - \frac{S^*}{S} \right) 
+ \varepsilon \beta I^* S^* \left( 2- \frac{S^*}{S} + \frac{I}{I^*} 
- \dfrac{I^*}{I}  \frac{L}{L^*} -  \dfrac{L^{*}}{L} \frac{I}{I^*} \frac{S}{S^*} \right)\\
&\quad + \varepsilon (1-k)\delta T^* \left( 1+ \frac{T}{T^*} - \dfrac{I^*}{I}  \frac{L}{L^*} 
- \dfrac{L^{*}}{L} \frac{T}{T^*} \right)+ b_1 k \delta T^* \left( \frac{T}{T^*} 
- \dfrac{I^*}{I} \frac{T}{T^*} \right)\\ 
&\quad - b_1 b_2 I  + b_1 I^* b_2 
+ m_3 \left(1-\dfrac{T^*}{T} \right)\left( \gamma I - b_3 T\right). 
\end{split}
\end{equation*}
After some simplifications, and using again relations 
$$
b_2 I^* = \varepsilon L^* + k \delta T^*,
\quad  
b_1 L^* = \beta I^* S^* + (1-k) \delta T^*, 
$$
we have 
\begin{equation*}
\begin{split}
_{t_{0}}^{C}D_{t}^{\alpha}V
&\leq \varepsilon  \mu S^* \left(2 - \frac{S}{S^*} - \frac{S^*}{S} \right) 
+ \varepsilon \beta I^* S^* \left( 3- \frac{S^*}{S} - \dfrac{I^*}{I}  \frac{L}{L^*} 
- \dfrac{L^{*}}{L} \frac{I}{I^*} \frac{S}{S^*} \right)\\
&\quad + \varepsilon (1-k)\delta T^* \left( 2 + \frac{T}{T^*} - \frac{I}{I^*} 
- \dfrac{I^*}{I}  \frac{L}{L^*} - \dfrac{L^{*}}{L} \frac{T}{T^*} \right)
+ b_1 k \delta T^* \left( 1 + \frac{T}{T^*} - \frac{I}{I^*} 
- \dfrac{I^*}{I} \frac{T}{T^*} \right)\\
&\quad + m_3 \gamma I - m_3 b_3 T - m_3 \dfrac{T^*}{T} \gamma I + m_3 T^* b_3. 
\end{split}
\end{equation*}
Let 
$$
m_3 = \displaystyle \frac{\delta T^* (b_1 k + \epsilon (1-k))}{\gamma I^*}.
$$ 
Using $b_3 T^* = \gamma I^*$, and simplifying the previous inequality, we have 
\begin{equation*}
\begin{split}
_{t_{0}}^{C}D_{t}^{\alpha}V &\leq \varepsilon  \mu S^* \left(2 
- \frac{S}{S^*} - \frac{S^*}{S} \right) + \varepsilon \beta I^* S^* 
\left( 3- \frac{S^*}{S} -  \dfrac{I^*}{I}  \frac{L}{L^*} 
- \dfrac{L^{*}}{L} \frac{I}{I^*} \frac{S}{S^*} \right)\\
&\quad + \varepsilon (1-k)\delta T^* \left( 3 - \frac{T^*}{T} \frac{I}{I^*} 
- \dfrac{I^*}{I}  \frac{L}{L^*} - \dfrac{L^{*}}{L} \frac{T}{T^*} \right)
+ b_1 k \delta T^* \left( 2 - \frac{T^*}{T} \frac{I}{I^*} 
- \dfrac{I^*}{I} \frac{T}{T^*} \right). 
\end{split}
\end{equation*}
Since the arithmetical mean is greater than or equal to the geometrical mean, 
then we have $_{t_{0}}^{C}D_{t}^{\alpha}V \leq 0$, with equality holding only
if $S^{*}=S$, $L^{*}= L$, $I^* = I$ and $T^* = T$. Therefore, by Theorem~\ref{uniform_stability}
of uniform asymptotic stability, we conclude that the endemic equilibrium $E^*$ \eqref{EE} 
is uniformly asymptotically stable in the interior of $\Omega$. 
\end{proof}


\section{Numerical simulations}
\label{sec:numsim}

In this section we study the dynamical behavior of our model \eqref{model_TB_fractional}, 
by variation of the noninteger order derivative $\alpha$. The parameter values used 
in the simulations can be found in Table~\ref{parameters}, with exception of $b=0.0005$, 
which we changed in order to get a value of $R_0$ bigger than one (endemic situation). 
Direct calculations, with these parameter values, give $R_0 =7.1343 > 1$ 
and the endemic equilibrium 
$$
E^{*}=(7779.3, 43512, 175.2, 78.4).
$$ 
We consider the initial conditions 
\begin{equation}
\label{init_cond}
S(0)=0.8,\quad L(0)=0.05,\quad I(0)=0.1,\quad T(0)=0.05
\end{equation}
and a fixed time step size of $h=2^{-8}$.
\begin{table}[h]
\caption{Parameter values taken from \cite{Yang201079}.}
\centering
\begin{tabular}{|c|l|c|} \hline 
Parameter & Description & Value\\ \hline 
$\Lambda$ & Recruitment rate  & 792.8571\\
$\beta$ & Transmission coefficient & $ 5 \times 10^{-5}$\\
$\mu$ & Natural death rate & $0.143$\\
$k$ & Treatment failure rate  & $0.15$ \\
$\delta$ & Rate at which treated individuals leave the $T$ compartment  & $1.5$\\   
$\epsilon$ & Rate at which latent individuals $L$ become infectious  & $0.00368$ \\
$\gamma$ & Treatment rate for infectious individuals $I$  & $0.7$\\
$\alpha_1$ & TB-induced death rate for infectious individuals $I$   & $0.3$ \\
$\alpha_2$ & TB-induced death rate for under treatment individuals $T$ & $0.05$ \\ \cline{1-3}
\end{tabular}
\label{parameters}
\end{table}

\begin{figure}
\centering
\includegraphics[width=0.48\textwidth]{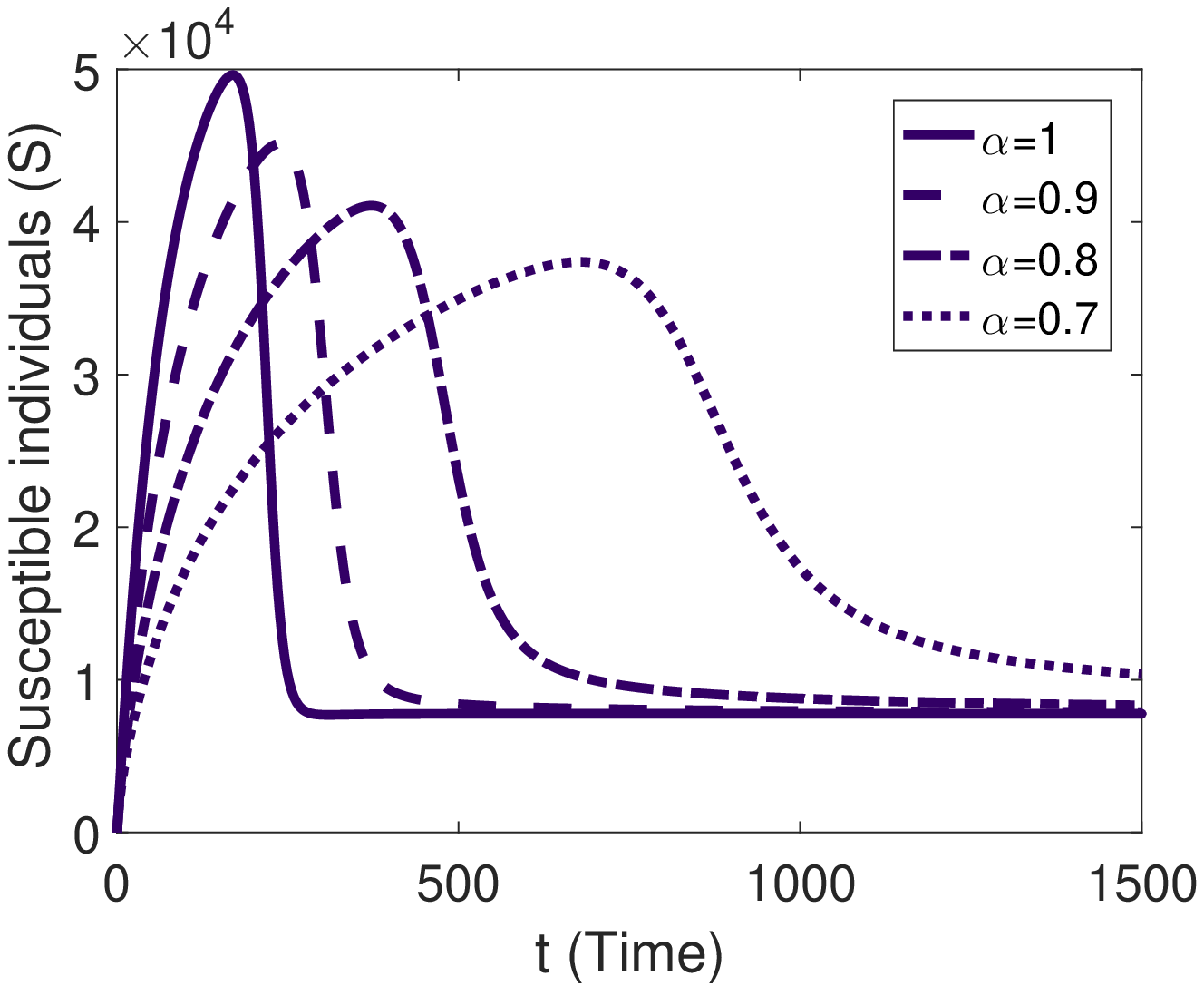}
\includegraphics[width=0.48\textwidth]{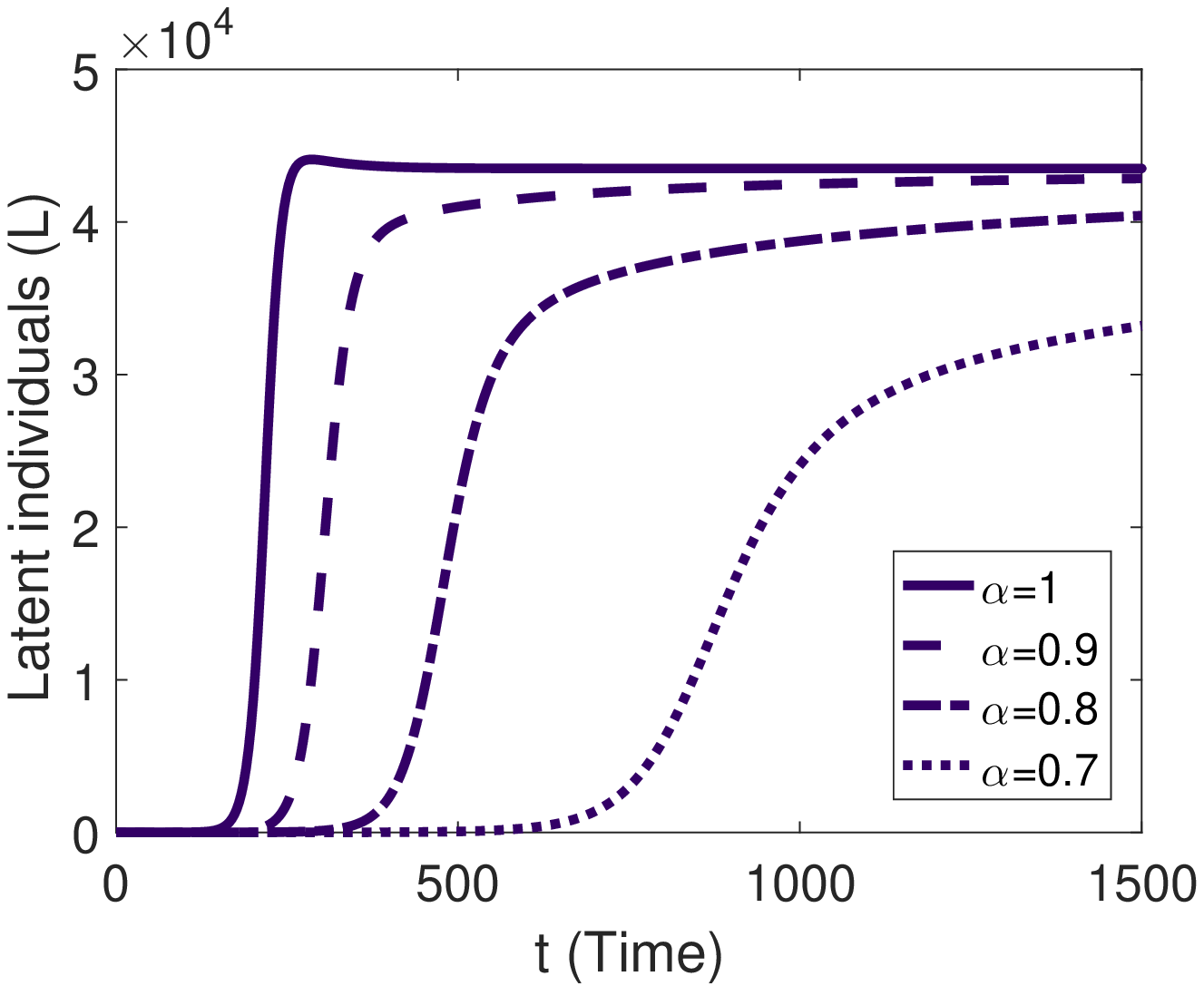}\\
\includegraphics[width=0.48\textwidth]{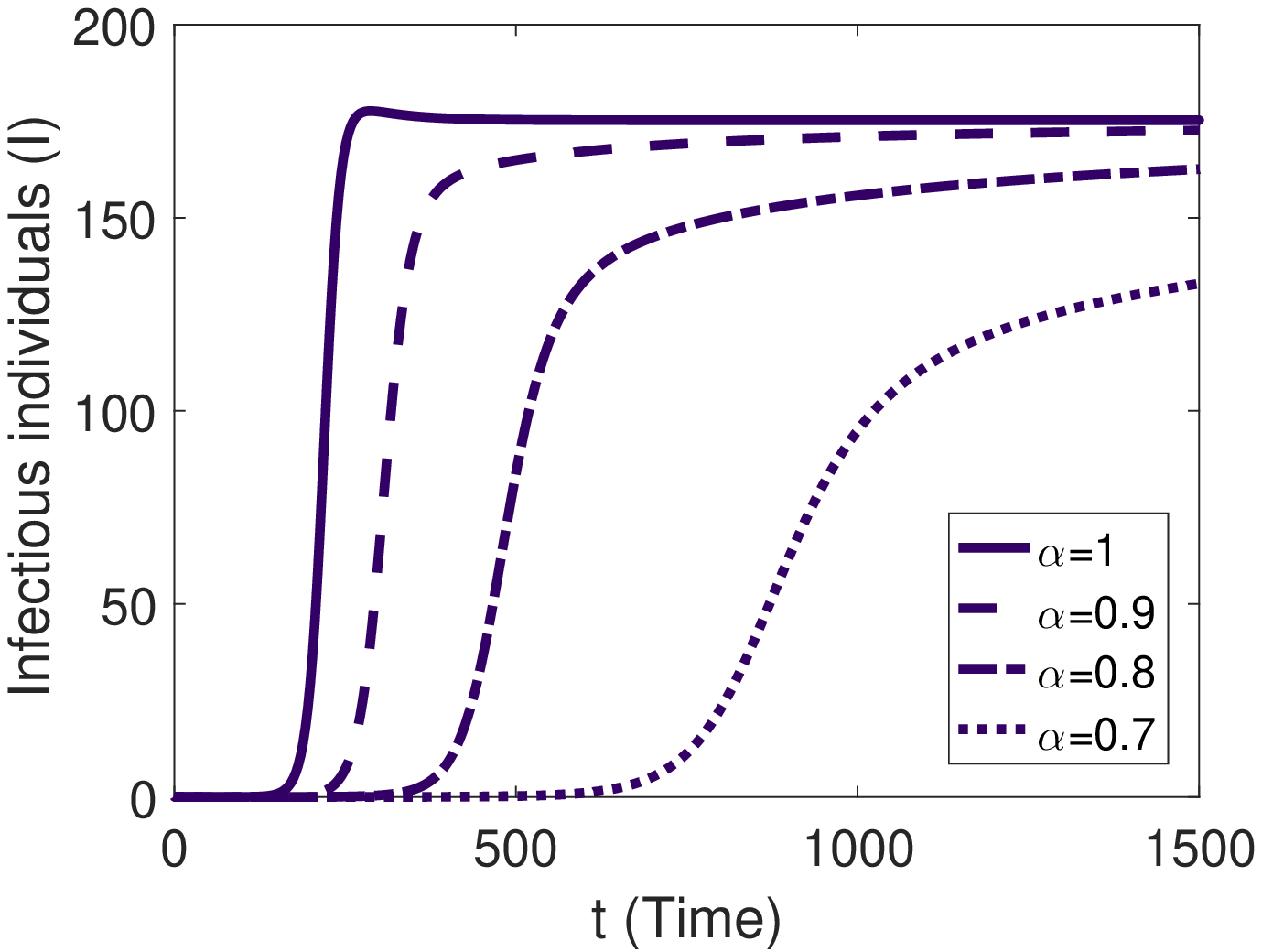}
\includegraphics[width=0.48\textwidth]{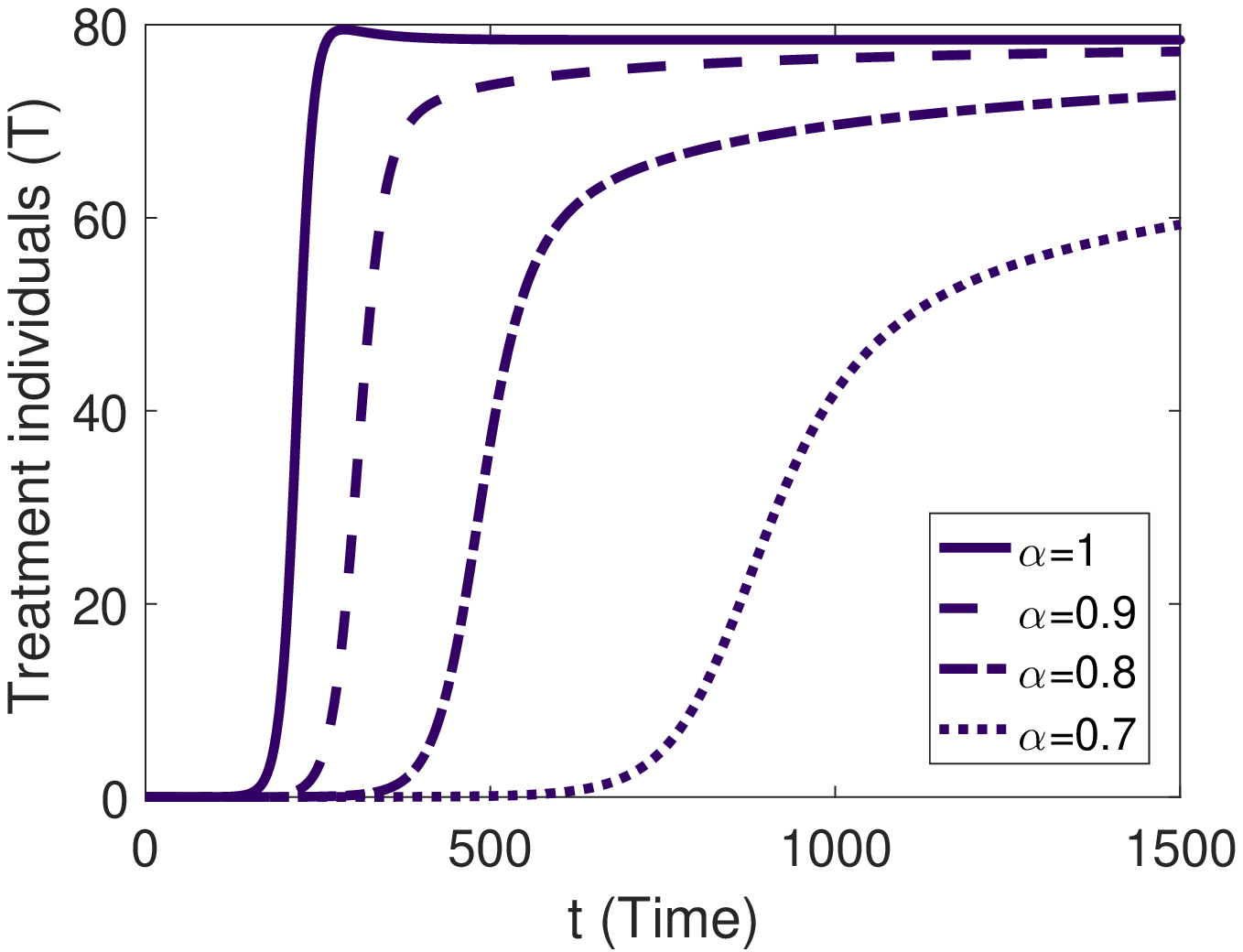}
\caption{Numerical simulation of the Caputo TB system 
\eqref{model_TB_fractional} with values as in Table~\ref{parameters}, 
with exception of $b=0.0005$, and initial conditions \eqref{init_cond}.}
\label{fig:simulation}
\end{figure}

\begin{figure}
\centering
\includegraphics[width=0.48\textwidth]{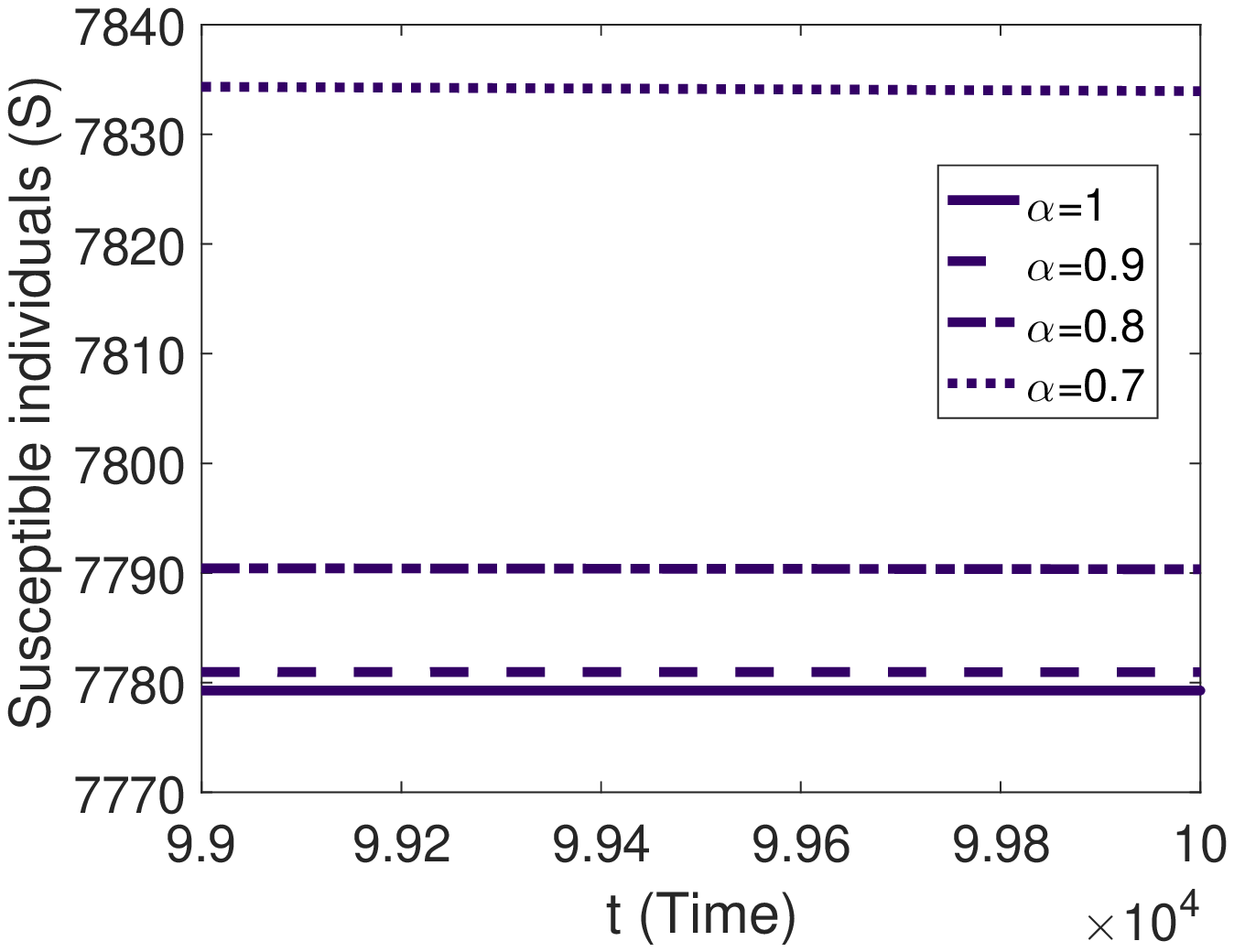}
\includegraphics[width=0.48\textwidth]{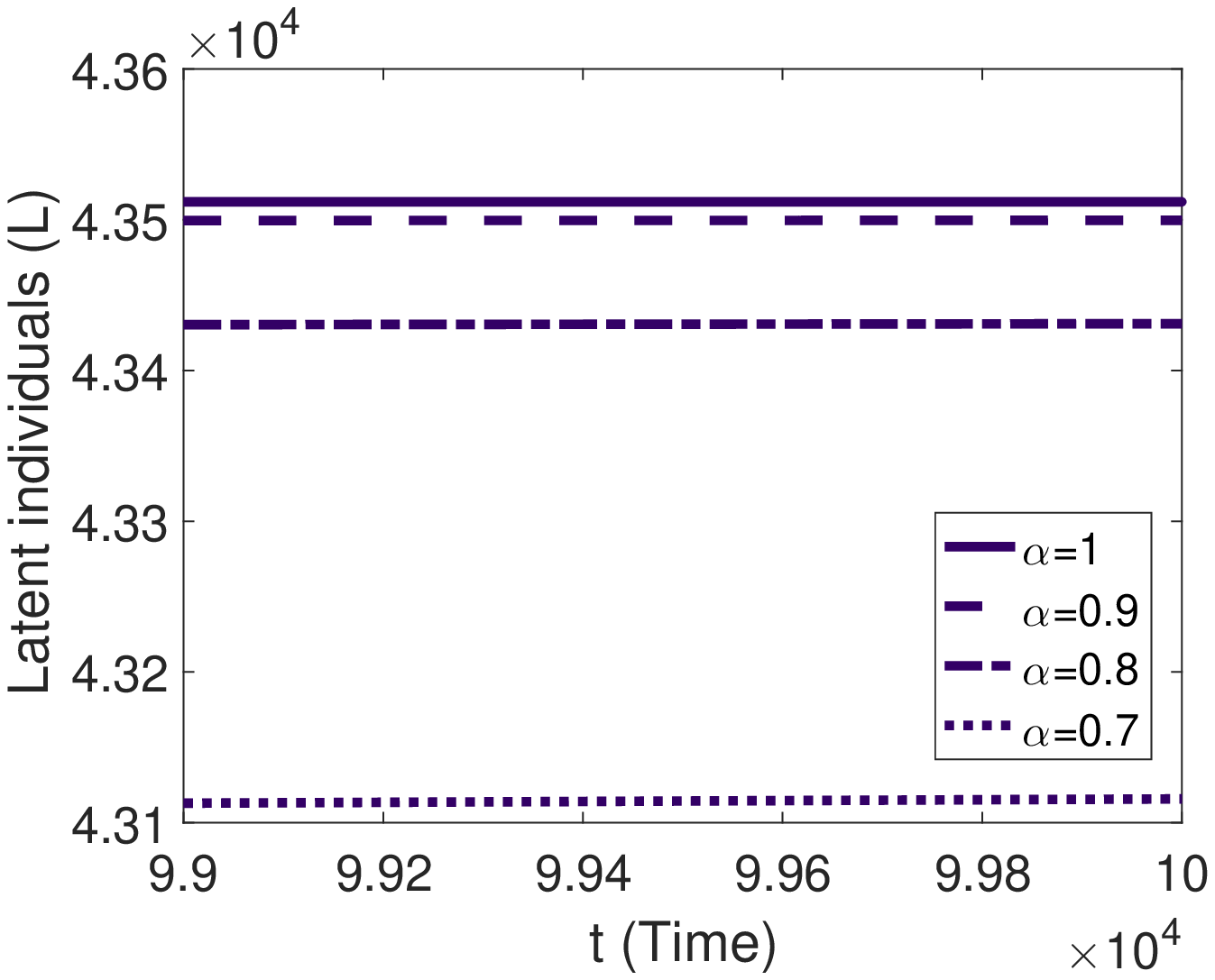}\\
\includegraphics[width=0.48\textwidth]{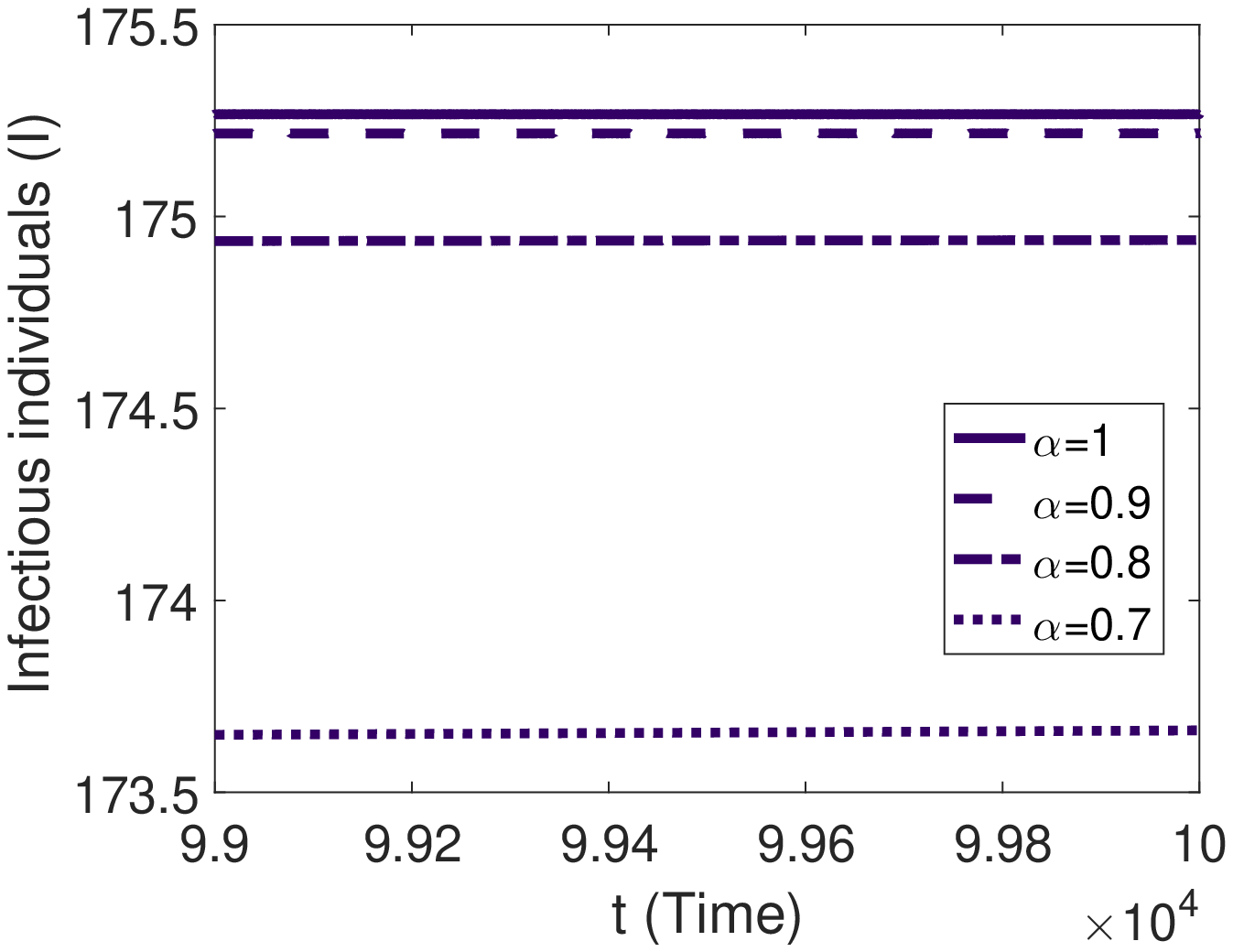}
\includegraphics[width=0.48\textwidth]{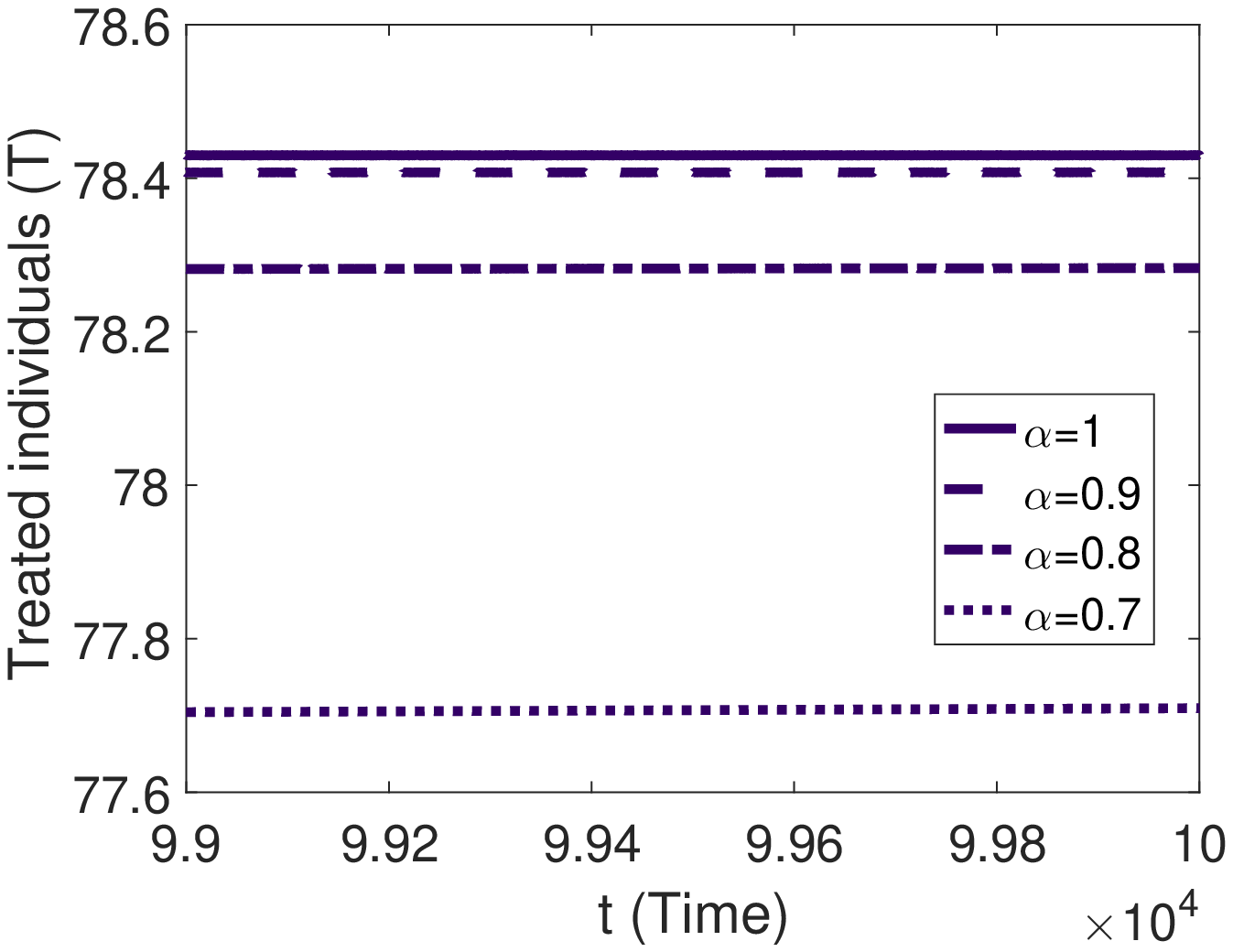}
\caption{Numerical simulation of the Caputo TB system  
with parameters and initial conditions as in Fig.~\ref{fig:simulation}, 
showing that the solutions converge to the endemic equilibrium 
$(7779.3, 43512, 175.2, 78.4)$ for $\alpha = 1$.}
\label{fig:simulation_long}
\end{figure}

For the numerical implementation of the fractional derivatives, 
we have used the Adams--Bashforth--Moulton scheme, which has been 
implemented in the \textsf{Matlab} code \textrm{fde12} 
by Garrappa \cite{Garrappa}. Regarding convergence and 
accuracy of the numerical method, we refer to \cite{[2]}.
The stability properties of the method implemented by \textrm{fde12} 
have been studied in \cite{[5]}.
We consider, without loss of generality, 
the fractional-order derivatives $\alpha = 1.0, 0.9, 0.8$ and $0.7$. 
In Figures~\ref{fig:simulation}--\ref{fig:simulation_long}, one can see 
that when $\alpha \rightarrow 1$ the solutions of our model converge 
to the solutions obtained in \cite{Yang201079}. The simulation results
confirm, numerically, the stability result of Theorem~\ref{thm:MR}.


\section{Conclusions}
\label{sec:conc}

In this work we proposed a Caputo fractional-order tuberculosis (TB) model. 
Existence of equilibrium points has been investigated and
the uniform asymptotic stability of the unique endemic equilibrium
proved via a suitable Lyapunov function. Our analytical results 
were complemented by numerical simulations in \textsf{Matlab},
illustrating the obtained stability result. 
The proposed fractional order model provides 
richer and more flexible results when compared with 
the corresponding integer-order TB model.


\begin{acknowledgement}
This research was partially supported by the
Portuguese Foundation for Science and Technology (FCT)
through the R\&D unit CIDMA, reference UID/MAT/04106/2013,
and by project PTDC/EEI-AUT/2933/2014 (TOCCATA),
funded by FEDER funds through COMPETE 2020 -- 
Programa Operacional Competitividade e
Internacionaliza\c{c}\~ao (POCI) 
and by national funds through FCT.
Silva is also grateful to the FCT post-doc 
fellowship SFRH/BPD/72061/2010; 
Wojtak to the FCT PhD fellowship PD/BD/128183/2016.
\end{acknowledgement}



\end{document}